\documentclass{article}

\title{Generalized Optimization: A First Step Towards Category Theoretic Learning Theory}

\author{
Dan Shiebler\\
daniel.shiebler@kellogg.ox.ac.uk\\
University of Oxford
}

\usepackage{tikz}

\usepackage[round]{natbib}
\usetikzlibrary{
  cd,
  math,
  decorations.markings,
  decorations.pathreplacing,
  positioning,
  arrows.meta,
  circuits.logic.US,
  shapes,
  calc,
  fit,
  trees,
  quotes}
\usetikzlibrary{decorations.pathmorphing}
\usetikzlibrary{decorations.markings}
\usetikzlibrary{decorations.pathreplacing}
\usetikzlibrary{arrows}
\usetikzlibrary{shapes.geometric}
\usepackage{amsthm}
\usepackage{amsmath}
\usepackage{amssymb}
\usepackage{hyperref}
\usepackage{bm}
\usepackage{bbm}
\usepackage{pdfpages}
\usepackage{float}
\usepackage{dsfont}
\usepackage{relsize}
\usepackage{todonotes}
\usepackage{enumitem}

\def\references#1{\vspace*{5mm}\noindent{References:}\list
{[\arabic{enumi}]}{\settowidth\labelwidth{[#1]}\leftmargin\labelwidth
\advance\leftmargin\labelsep
\usecounter{enumi}}
\def\newblock{\hskip .11em plus .33em minus .07em}
\sloppy\clubpenalty4000\widowpenalty4000
\sfcode`\.=1000\relax}

\newcommand{\optimizer}{\mathbf{Optimizer}}

\newcommand{\db}{\mathbf{D}}
\newcommand{\cb}{\mathbf{C}}
\newcommand{\rl}{\mathbb{R}}

\newcommand{\euc}{\mathbf{Euc}}

\newcommand{\base}{\mathbf{Base}}
\newcommand{\basedomain}{(\base, X)}

\newcommand{\objective}{\mathbf{Objective}}
\newcommand{\objectivei}{\objective_{I}}
\newcommand{\objectiveo}{\objective_{\perp}}

\newcommand{\ring}{r}
\newcommand{\polyr}{\mathbf{Poly}_{\ring}}
\newcommand{\polyz}{\mathbf{Poly}_{\mathbb{Z}}}

\newtheorem{proposition}{Proposition}

\newtheorem{definition}{Definition}[section]

\begin{document}
\maketitle

\begin{abstract}
The Cartesian reverse derivative is a categorical generalization of reverse-mode automatic differentiation. We use this operator to generalize several optimization algorithms, including a straightforward generalization of gradient descent and a novel generalization of Newton's method. We then explore which properties of these algorithms are preserved in this generalized setting. First, we show that the transformation invariances of these algorithms are preserved: while generalized Newton's method is invariant to all invertible linear transformations, generalized gradient descent is invariant only to orthogonal linear transformations. Next, we show that we can express the change in loss of generalized gradient descent with an inner product-like expression, thereby generalizing the non-increasing and convergence properties of the gradient descent optimization flow.
Finally, we include several numerical experiments to illustrate the ideas in the paper
and demonstrate how we can use them to optimize polynomial functions over an ordered ring.
\end{abstract}

\section{Background}\label{introduction}
Given a convex differentiable function $l: \rl^{n} \rightarrow \rl$, there are many algorithms that we can use to minimize it. For example, if we pick a step size $\alpha$ and a starting point $x_0 \in \rl^n$ we can apply the gradient descent
algorithm in which we repeatedly iterate $x_{t+1} = x_{t} - \alpha * \nabla l(x_t)$. For small enough $\alpha$ this strategy is guaranteed to get close to the $x$ that minimizes $l$ \citep{convexopt}. 

Algorithms like gradient descent are often useful even when $l$ is non-convex. For example, under relatively mild conditions we can show that taking small enough gradient descent steps will never increase the 
value of any differentiable $l: \rl^n \rightarrow \rl$ \citep{convexopt}. The modern field of deep learning consists largely of applying gradient descent and other algorithms that can be efficiently computed with reverse-mode automatic differentiation to optimize non-convex functions \citep{lecun2015deep}.

While gradient descent is particularly popular, it is not the only gradient-based optimization algorithm that is widely used in practice. Both the momentum and Adagrad algorithms use placeholder variables that store information from previous gradient updates to improve stability. Newton's method, which rescales the gradient with the inverse Hessian matrix, is popular for many applications but less commonly used in deep learning due to the difficulty of efficiently implementing it with reverse-mode automatic differentiation. Each of these algorithms have different invariance properties that affect their stability under dataset transformations: for example, Newton's method enjoys an invariance to linear rescaling that gradient descent lacks.

Given the utility of these algorithms it is natural to explore when they can be generalized beyond differentiable functions. For example, if we replace the gradient in gradient descent with a representative of the subgradient, a simple generalization of the gradient for non-differentiable functions $l: \rl^{n} \rightarrow \rl$, the convergence and stability properties of gradient descent are preserved \citep{convexopt}. Going one step further, some authors have begun to explore generalizations of differentiation beyond Euclidean spaces. \cite{cockett2019reverse} introduce Cartesian reverse derivative categories in which we can define an operator that shares certain properties with reverse-mode automatic differentiation (RD.1 to RD.7 in Definition 13 of \cite{cockett2019reverse}). Reverse derivative categories are remarkably general: the category of Euclidean spaces and differentiable functions between them is of course a reverse derivative category, as are the categories of polynomials over semirings and the category of Boolean circuits \citep{Wilson_2021}.

In this paper we explore whether the convergence and invariance properties of optimization algorithms built on top of the gradient and Hessian generalize to optimization algorithms built on top of \cite{cockett2019reverse}'s Cartesian reverse derivative. Our contributions are as follows:
\begin{itemize}
    %
    \item We use \cite{cockett2019reverse}'s Cartesian reverse derivative to define generalized analogs of several optimization algorithms, including a novel generalization of Newton's method.
    \item We derive novel results on the transformation invariances of these generalized algorithms. 
    \item We define the notion of an optimization domain over which we can apply these generalized algorithms and study their convergence properties. 
    \item We characterize the properties that an optimization domain must satisfy in order to support generalized gradient-based optimization and we provide novel results that the optimization domain of polynomials over ordered rings satisfies these properties.
    \item We include several numerical experiments to illustrate the ideas in the paper and demonstrate how we can use them to optimize polynomial functions over an ordered ring. The code to run these experiments is on Github at tinyurl.com/ku3pjz56.
\end{itemize}

\section{Standard Optimization}\label{standardoptimization}

As we described in Section \ref{introduction}, gradient descent optimizes an objective function $l: \rl^n \rightarrow \rl$ by starting at a point $x_0 \in \rl^n$ and progressing the discrete dynamical system $x_{t+1} = x_{t} - \alpha * \nabla l(x_t)$. Rewriting this as $x_{t+\alpha} = x_{t} - \alpha * \nabla l(x_t)$ and taking the $lim_{\alpha \rightarrow 0}$ of this system yields the differential equation $\frac{\partial x}{\partial t}(t) = -\nabla l(x(t))$, which we can think of as the continuous limit of gradient descent.
More generally we have:
\begin{definition}\label{optimizerdefinition}
An \textbf{optimizer} for $l: \rl^n \rightarrow \rl$ with \textbf{dimension} $k$ is a continuous function $d: \rl^{kn} \rightarrow \rl^{kn}$.
\end{definition}
Intuitively, an optimizer defines both a continuous and a discrete dynamical system
\begin{itemize}
    \item \textbf{Continuous System}: $\left(\frac{\partial x}{\partial t}(t), \frac{\partial y}{\partial t}(t), \cdots\right) = d(x(t), y(t), \cdots)$
    \item \textbf{Discrete System}: $\left(x_{t + 1},y_{t + 1}, \cdots) = (x_{t },y_t, \cdots) + \alpha d(x_t, y_t, \cdots\right)$
\end{itemize}
Note that the discrete dynamical system is the Euler's method discretization of the continuous system. We can think of an optimizer with dimension $k > 1$ as using information beyond the previous value $x_t$ to determine $x_{t+1}$.

In practice we usually work with optimizers that define dynamical systems in which $l(x(t))$ and $l(x_t)$ get closer to the minimum value of $l$ as $t$ increases. Given $l: \rl^{n}\rightarrow \rl$ we can construct the following examples:
\begin{itemize}
    \item 
    The \textbf{gradient descent optimizer} for $l: \rl^n \rightarrow \rl$ is $d(x) = -\nabla l(x)$. This optimizer has dimension $1$.
    \item
    The \textbf{Newton's method optimizer} for $l: \rl^n \rightarrow \rl$ is $d(x) = -\nabla^2 (l (x))^{-1}\nabla l(x)$. This optimizer has dimension $1$ and uses the inverse Hessian matrix to increase the stability of each update.
    \item 
    The \textbf{momentum optimizer}
    for $l: \rl^n \rightarrow \rl$ is $d(x,y) = (y, -y - \nabla l(x))$. This optimizer has dimension $2$ and uses a placeholder variable to track the value of the previous update steps and simulate the momentum of a ball rolling down a hill. 
    \item
    The \textbf{Adagrad optimizer}
    for $l: \rl^n \rightarrow \rl$ is $d(x,y) = (-\nabla l(x) / \sqrt{y}, \nabla l(x)^2)$. This optimizer has dimension $2$ and uses a placeholder variable to reweight updates based on the magnitude of previous updates \citep{duchi2011adaptive}. 
\end{itemize}


%

\subsection{Optimization Schemes}\label{continuousoptimizationalgorithm}
\begin{definition}\label{optimizationschemedefinition}
An \textbf{optimization scheme}
$u: (\rl^{n} \rightarrow \rl) \rightarrow (\rl^{kn} \rightarrow \rl^{kn})$
 is a family of functions (indexed by $n$) that maps objectives $l: \rl^{n} \rightarrow \rl$ to optimizers
$d: \rl^{kn} \rightarrow \rl^{kn}$.
\end{definition}
For example, the gradient descent optimization scheme is $u(l)(x) = -\nabla l (x)$ and the momentum optimization scheme is $u(l)(x, y) = (y, -y - \nabla l(x))$.

In some situations we may be able to improve the convergence rate of the dynamical systems defined by optimization schemes by precomposing an invertible function $f: \rl^m \rightarrow \rl^n$. That is, rather than optimize the function $l: \rl^{n} \rightarrow \rl$ we optimize $l \circ f: \rl^{m} \rightarrow \rl$. However, for many optimization schemes there are classes of transformations to which they are invariant: applying any such transformation to the data cannot change the trajectory.
\begin{definition}
Suppose $f: \rl^{m} \rightarrow \rl^n$ is an invertible transformation and write $f_{k}$ for the map $(f \times f \times \cdots): \rl^{km} \rightarrow \rl^{kn}$. The optimization scheme $u$ is \textbf{invariant} to $f$ if $u(l \circ f) = f_{k}^{-1} \circ u(l) \circ f_{k}$.
\end{definition}
\begin{proposition}\label{invarianceprop}
Recall that an invertible linear transformation is a function $f(x) = Ax$ where the matrix $A$ has an inverse $A^{-1}$ and an orthogonal linear transformation is an invertible linear transformation where $A^{-1} = A^{T}$. Newton's method is invariant to all invertible linear transformations, whereas both gradient descent and momentum  are invariant to orthogonal linear transformations.
(Proof in Appendix \ref{invarianceproof}).
\end{proposition}
Note that Adagrad is not invariant to orthogonal linear transformations due to the fact that it tracks a nonlinear function of past gradients (sum of squares). In order to interpret these invariance properties it is helpful to consider how they affect the discrete dynamical system defined by an optimization scheme.
\begin{proposition}\label{eulersmethodprop}
Given an objective function $l: \rl^n \rightarrow \rl$ and an optimization scheme $u: (\rl^{n} \rightarrow \rl) \rightarrow (\rl^{kn} \rightarrow \rl^{kn})$ that is invariant to the invertible linear  function $f: \rl^m \rightarrow \rl^n$, the discrete system defined by the optimizer $u(l \circ f)$:
\begin{gather*}
y_{t + 1} = y_t + \alpha u(l \circ f)(y_t)
\end{gather*}
cannot converge faster than the discrete system defined by the optimizer $u(l)$:
\begin{gather*}
x_{t + 1} = x_t + \alpha u(l)(x_t)
\end{gather*}
(Proof in Appendix \ref{eulersmethodproof}).
\end{proposition}

Propositions 1 and 2 together give some insight into why Newton's method can perform so much better than gradient descent for applications where both methods are computationally feasible \citep{convexopt}. Whereas gradient descent can be led astray by bad data scaling, Newton's method steps are always scaled optimally and therefore cannot be improved by data rescaling. 

It is important to note that Proposition \ref{eulersmethodprop} only applies to linear transformation functions $f$. Since Euler's method is itself a linear method, it does not necessarily preserve non-linear invariance properties.

\section{Generalized Optimization}\label{generalized}

In this section we assume readers have a basic familiarity with Category Theory. We recommend that readers who would like a detailed introduction to the field check out ``Basic Category Theory'' \citep{leinster2016basic} or ``Seven Sketches in Compositionality'' \citep{fong2018seven}. 

\subsection{Cartesian and Reverse Derivative Categories}\label{differential}

We begin by recalling some basic terminology and notation. A \textbf{category} is a collection of objects and morphisms between them. Morphisms are closed under an associative composition operation, and each object is equipped with an identity morphism. An example category is the collection of sets (objects) and functions (morphisms) between them. We call the set of morphisms between the objects $A$ and $B$ in the category $\cb$ the \textbf{hom-set} $\cb[A, B]$. 

A category is \textbf{Cartesian} when there exists a product operation $\times$ that allows us to combine objects, a pairing operation $\langle -, - \rangle$ that allows us to combine morphisms, projection maps $\pi_0: A \times B \rightarrow A, \pi_1: A \times B \rightarrow B$ and a terminal object $*$ such that every object $A$ is equipped with a unique map $!_{A}: A \rightarrow *$ from $A$ to the terminal object $*$. Given an object $A$ or morphism $f$ in a Cartesian category, in this work we write $A^{k}$ and $f^{k}$ to respectively denote $A$ and $f$ tensored with themselves $k$ times.

Now recall the following definition from  \cite{blute2009cartesian, cockett2019reverse}:
\begin{definition}
A \textbf{Cartesian left additive category} is a Cartesian category $\cb$ in which the hom-set of each pair of objects $B, C$ is a commutative monoid, with addition operation $+$ and zero maps (additive identities) $0_{BC}: B \rightarrow C$, such that:
\begin{itemize}
    \item For any morphism $h: A \rightarrow B$ and $f,g: B\rightarrow C$ we have:
    \begin{gather*}
        (f + g) \circ h = (f \circ h) + (g \circ h): A \rightarrow C
        \qquad
        0_{BC} \circ h= 0_{AC}: A\rightarrow C
    \end{gather*}
    \item For any projection map $\pi_i: C \rightarrow D$ and $f,g: B\rightarrow C$ we have:
    \begin{gather*}
        \pi_i \circ (f + g) = (\pi_i \circ f) + (\pi_i \circ g): B \rightarrow D
        \qquad
        \pi_i \circ 0_{BC} = 0_{BD}: B \rightarrow D
    \end{gather*}
\end{itemize}
We write $0_A$ for the additive identity of the hom-set $\cb[*, A]$.
\end{definition}
%
Intuitively, in a Cartesian left additive category we can add morphisms in a way that is compatible with postcomposition and the Cartesian structure. Certain Cartesian left additive categories are equipped with additional structure that behaves similarly to derivatives:
\begin{definition}
A \textbf{Cartesian reverse derivative category} is a Cartesian left-additive category $\cb$ equipped with a \textbf{Cartesian reverse derivative combinator} $R$ that assigns to each morphism $f : A \rightarrow B$ in $\cb$ a morphism $R[f] : A \times B \rightarrow A$ in $\cb$ such that $R$ satisfies the following equations (Definition 13 of \citep{cockett2019reverse}):
\begin{itemize}[labelindent=4em,labelsep=0.5cm,leftmargin=*]
    \item[\textbf{RD.1}] $R[f + g] = R[f] + R[g]$ and $R[0] = 0$;
    \item[\textbf{RD.2}] $R[f] \circ \langle a, b + c \rangle = R[f] \circ \langle a, b \rangle + R[f] \circ \langle a, c \rangle$
    and
    $R[f] \circ \langle a, 0 \rangle = 0$; 
    \item[\textbf{RD.3}] $R[1] = \pi_1, R[\pi_0] = \iota_0 \circ \pi_1$ and $R[\pi_1] = \iota_1 \circ \pi_1$;
    \item[\textbf{RD.4}] $R[\langle f, g \rangle] = R[f] \circ (1 \times \pi_0) + R[g] \circ (1 \times \pi_1)$ and $R[!_A] = 0$;
    \item[\textbf{RD.5}] $R[g \circ f] = R[f] \circ (1 \times R[g]) \circ \langle \pi_0, \langle f \circ \pi_0, \pi_1 \rangle \rangle$;
    \item[\textbf{RD.6}] $
    \pi_1
    \circ
    R[R[R[f]]]
    \circ
    (\iota_0 \times 1)
    \circ
    \langle 1 \times \pi_0,
    0 \times \pi_1 \rangle = 
    R[f] \circ (1 \times \pi_1)$;
    \item[\textbf{RD.7}] \begin{gather*}
        \pi_1 \circ
    R[R[\pi_1 \circ R[R[f]] \circ (\iota_0 \times 1)]
    ] \circ (\iota_0 \times 1)
    =\\
    \langle \pi_0 \times \pi_0, \pi_1 \times \pi_1 \rangle 
    \circ
    \pi_1 \circ
    R[R[\pi_1 \circ R[R[f]] \circ (\iota_0 \times 1)]
    ] \circ (\iota_0 \times 1)
    \end{gather*}
\end{itemize}
where $\iota_0: A \rightarrow A \times B$ and $\iota_1: B \rightarrow A \times B$ are the Cartesian injection maps \citep{cockett2019reverse}.  
\end{definition}
The conditions \textbf{RD.1} to \textbf{RD.7} mirror the properties of the derivative operation. For example, $R$ must commute with addition (\textbf{RD.1}) and compose according to a chain rule (\textbf{RD.5}).
\begin{definition}
A \textbf{Cartesian differential category} $\cb$ is a Cartesian left-additive category equipped with a \textbf{Cartesian derivative combinator} $D$ that assigns to each morphism $f : A \rightarrow B$ in $\cb$ a morphism $D[f] : A \times A \rightarrow B$ in $\cb$ such that $D$ satisfies the following equations (Definition 4 in \cite{cockett2019reverse}, adapted from Definition 2.1.1 in   \citep{blute2009cartesian}):
\begin{itemize}[labelindent=4em,labelsep=0.5cm,leftmargin=*]
    \item[\textbf{CDC.1}] $D[f + g] = D[f] + D[g]$ and $D[0] = 0$;
    \item[\textbf{CDC.2}] $D[f] \circ \langle a, b + c\rangle = D[f] \circ \langle a, b\rangle + D[f] \circ \langle a, c \rangle$
    and
    $D[f] \circ \langle a, 0 \rangle  = 0$;
    \item[\textbf{CDC.3}] $D[1] = \pi_1, D[\pi_0] = \pi_0 \circ \pi_1$ and $D[\pi_1] = \pi_1 \circ \pi_1$;
    \item[\textbf{CDC.4}] $D[\langle f, g \rangle] = \langle D[f], D[g]\rangle$;
    \item[\textbf{CDC.5}] $D[g\circ f] = D[g] \circ \langle f \circ \pi_0, D[f]\rangle$;
    \item[\textbf{CDC.6}] $D[D[f]] \circ \langle\langle a,b\rangle,\langle 0,c \rangle \rangle = D[f] \circ \langle a, c\rangle$;
    \item[\textbf{CDC.7}] $D[D[f]] \circ \langle \langle a,b \rangle, \langle c,d\rangle \rangle = D[D[f]] \circ \langle \langle a,c \rangle, \langle b,d \rangle \rangle$.
\end{itemize}
\end{definition}

By Theorem 16 in \citep{cockett2019reverse}, every Cartesian reverse derivative category $\cb$ is also a Cartesian differential category where for any morphism $f: A \rightarrow B$ in $\cb$:
\begin{gather*}
    D[f] = \pi_1 \circ R[R[f]] \circ (\langle id_A, 0_{AB}\rangle \times id_A): A \times A \rightarrow B
\end{gather*}
Going forward, when we refer to the Cartesian derivative combinator $D$ of a Cartesian reverse derivative category this is the construction that we are referring to. 

The canonical example of a Cartesian reverse derivative category that we will consider is the category $\euc$ of Euclidean spaces and infinitely differentiable maps between them. The terminal object $*$ in $\euc$ is $\rl^0$, the Cartesian reverse derivative of the map $f: \rl^a \rightarrow \rl^b$ is $R[f](x, x') = J_f(x)^{T} x'$, and the Cartesian derivative of $f$ is $D[f](x, x') = J_f(x)  x'$ where $J_f(x)$ is the Jacobian of $f$ evaluated at $x \in \rl^a$. Recall that the Jacobian of $f: \rl^a \rightarrow \rl^b$ is a $b \times a$ matrix whose $i,j$th element is $\frac{\partial f_i}{\partial x_j}$.

%
%
As another example, given a commutative ring $\ring$ we can form the category $\polyr$ in which objects are natural numbers and the morphisms from $n$ to $m$ are tuples of $m$ polynomials with $n$ variables and coefficients in $\ring$. That is, a morphism $P: n \rightarrow m$ is a map
$P(x) = \left( p_1(x), \cdots, p_m(x) \right)$
where $x = (x_1, \cdots, x_n)$ and $p_i$ is a polynomial. $\polyr$ is a Cartesian reverse derivative category in which the terminal object $*$ is $0$ and the reverse derivative of 
$P(x) = \left( p_1(x), \cdots, p_m(x) \right)$ is:
\begin{gather*}
R[P](x, x') = \left(
\sum_{i=1}^m \frac{\partial p_i}{\partial x_1}(x) x'_i,
\cdots,
\sum_{i=1}^m \frac{\partial p_i}{\partial x_n}(x) x'_i
\right)
\end{gather*}
%
%
where $\frac{\partial p_i}{\partial x_j}(x)$ is the formal derivative of $p_i$ in $x_j$, evaluated at $x$ \citep{cockett2019reverse}.


The \textbf{linear maps} in a Cartesian reverse derivative category $\cb$ are those for which $D[f] = f \circ \pi_1$. The linear maps of $\cb$ form a subcategory of $\cb$ equipped with a stationary on objects involution $(\textunderscore)^\dagger: \cb^{op} \rightarrow \cb$ such that for any linear map $f$ we have $R[f] = f^{\dagger} \circ \pi_1$ \citep{cockett2019reverse}.

The linear maps in $\euc$ are exactly the linear maps in the traditional sense, and given a linear map $f: \rl^a \rightarrow \rl^b$ in $\euc$ where $f(x) = Mx$ we have $f^{\dagger}: \rl^b \rightarrow \rl^a$ where $f^{\dagger}(x) = M^{T}x$. That is, $\dagger$ is a generalization of the transpose. Similarly, the linear maps in $\polyr$ are those that can be expressed as $P(x) = \left( p_1(x), \cdots, p_m(x) \right)$ where $p_i(x) = \sum_{i=1}^n r_i x_i$ for $r_i \in r$ \citep{cockett2019reverse}.


\subsection{Optimization Domain}\label{optimizationdomain}
%
%


%
\begin{definition}
%
An \textbf{optimization domain} is a tuple $\basedomain$ such that each morphism $f: A \rightarrow B$ in the Cartesian reverse derivative category $\base$ has an additive inverse $-f$ and each homset $\cb[*,A]$ out of the terminal object $*$ is further equipped with a multiplication operation $fg$ and a multiplicative identity map $1_{A}: * \rightarrow A$ to form a commutative ring with the left additive structure $+$.
$X$ is an object in  $\base$ such that the homset $f \in \cb[*,X]$ is further equipped with a total order $f \leq g$ to form an ordered commutative ring.
\end{definition}

Given an optimization domain $\basedomain$ the object $X$ represents the space of objective values to optimize and we refer to morphisms into $X$ as \textbf{objectives}. We abbreviate the map $1_B \circ !_{A}: A \rightarrow B$ as $1_{AB}$, where $!_{A}: A \rightarrow *$ is the unique map into the terminal object $*$.
For example, the objectives in the \textbf{standard domain} $(\euc, \rl)$ are functions $l: \rl^n \rightarrow \rl$. If $r$ is an ordered commutative ring then we can form the \textbf{$r$-polynomial domain} $(\polyr, 1)$ in which objectives are $r$-polynomials $l_P: n \rightarrow 1$.

%
\begin{definition}\label{boundedbelow}
We say that an objective $l: A \rightarrow X$ is \textbf{bounded below} in $\basedomain$ if there exists some $x: * \rightarrow X$ such that for any $a: * \rightarrow A$ we have $x \leq l \circ a$.
\end{definition}
In both the standard and $r$-polynomial domains an objective is bounded below if its image has an infimum.

\subsubsection{Generalized Gradient and Generalized $n$-Derivative}\label{generalizedgradientderivativesection}

\begin{definition}
The \textbf{generalized gradient} of the objective $l: A \rightarrow X$ in $\basedomain$ is $R[l]_1: A \rightarrow A$ where:
\begin{gather*}
    R[l]_1 = R[l] \circ \langle id_A, 1_{AX} \rangle
\end{gather*}
\end{definition}
In the standard domain the generalized gradient of $l: \rl^n \rightarrow \rl$ is just the gradient $R[l]_1(x) = \nabla l(x)$ and in the $r$-polynomial domain the generalized gradient of $l_P: n \rightarrow 1$  is
$
R[l_P]_1(x) = \left(\frac{\partial l_P}{\partial x_1}(x), \cdots, \frac{\partial l_P}{\partial x_n}(x) \right)
$
where $\frac{\partial l_P}{\partial x_i}$ is the formal derivative of the polynomial $l_P$ in $x_i$.
\begin{definition}
The \textbf{generalized $n$-derivative} of the morphism $f: X \rightarrow A$ in $\basedomain$ is $D_n[f]: X \rightarrow A$ where:
\begin{gather*}
    D_1[f] = D[f] \circ \langle id_X, 1_{XX}\rangle
    \qquad
    D_n[f] = D[D_{n-1}[f]] \circ\langle id_X, 1_{XX}\rangle
\end{gather*}
\end{definition}
In the standard domain the generalized $n$-derivative of $f: \rl \rightarrow \rl$ is the $n$-derivative $f^{(n)} = \frac{\partial^{n} f}{\partial x^{n}}$ and in the $r$-polynomial domain the generalized $n$-derivative of $l_P: 1 \rightarrow 1$ is the formal $n$-derivative $\frac{\partial^{n} l_P}{\partial x^{n}}$.

The derivative over the reals has a natural interpretation as a rate of change. We can generalize this as follows:
\begin{definition}\label{nsmooth}
We say that a morphism $f: X \rightarrow X$ in $\base$ is \textbf{$n$-smooth} in $\basedomain$ if whenever $D_{k}[f] \circ t \geq 0_X: * \rightarrow X$ for all
$t_1 \leq t \leq  t_2: * \rightarrow X$ and $k\leq n$ we have that $f \circ t_1 \leq f \circ t_2: * \rightarrow X$. 
\end{definition}
Intuitively, $f$ is $n$-smooth if it cannot decrease on any interval over which its generalized derivatives of order $n$ and below are non-negative. Some examples include:
\begin{itemize}
    \item 
    Any map $f: \rl \rightarrow \rl$ is $1$-smooth in the standard domain by the mean value theorem.
    \item  
    %
    When $r$ is a dense subring of a real-closed field then any polynomial $l_P: 1 \rightarrow 1$ is $1$-smooth in the $r$-polynomial domain \citep{stackexchangering}.
    \item  
    %
    For any $r$, the polynomial $l_P = \sum_{k=0}^{n} c_k t^k : 1 \rightarrow 1$ of degree $n$ is $n$-smooth in the $r$-polynomial domain since for any $t_1$ we can use the binomial theorem to write:
    \begin{gather*}
        l_P(t) =
        %
        \sum_{k=0}^{n} c_k t^k =
        %
        \sum_{k=0}^{n} c_k (t_1 + (t - t_1))^k =
        %
        l_P(t_1) + \sum_{k=1}^{n} c'_k (t - t_1)^k
    \end{gather*}
    where $c'_k$ is a constant such that $(c'_k)(k!) = D_{k}[l_P](t_1)$. Note that $c'_k$ must exist by the definition of the formal derivative of $l_P$, and must be non-negative if $D_{k}[l_P](t_1)$ is non-negative.
\end{itemize}


\subsection{Optimization Functors}\label{continuousoptimizationfunctorsection}

In this section we generalize optimization schemes (Section \ref{continuousoptimizationalgorithm}) to arbitrary optimization domains. This will enable us to characterize the invariance properties of our generalized optimization schemes in terms of the categories out of which they are functorial. Given an optimization domain $\basedomain$ we can define the following categories:
\begin{definition}
The objects in the category $\objective$ over the optimization domain $\basedomain$ are objectives $l: A \rightarrow X$ such that there exists an inverse function $R[l]_1^{-1}: A \rightarrow A$ where $R[l]_1^{-1} \circ R[l]_1 = R[l]_1 \circ R[l]_1^{-1} = id_A: A \rightarrow A$, and the morphisms between $l: A \rightarrow X$ and $l': B \rightarrow X$ are morphisms $f: A \rightarrow B$ where $l' \circ f = l$.
\end{definition}
\noindent Note that $\objective$ is a subcategory of the slice category $\base / X$.

In the standard domain the objects in $\objective$ are objectives $l: \rl^n \rightarrow \rl$ such that the function  $\nabla l: \rl^n \rightarrow \rl^n$ is invertible. In the $r$-polynomial domain, the objects in $\objective$ are $r$-polynomials $l_P: n \rightarrow 1$ such that the function $\langle \frac{\partial l_P}{\partial x_1}, \cdots, \frac{\partial l_P}{\partial x_n} \rangle: n \rightarrow n$ is invertible.
\begin{definition}\label{optbasedefinition}
%
A \textbf{generalized optimizer} over the optimization domain $\basedomain$ with \textbf{state space} $A \in \base$ and \textbf{dimension} $k \in \mathbb{N}$ is an endomorphism $d: A^{k}\rightarrow A^{k}$ in $\base$. The objects in the category $\optimizer$ over $\basedomain$ are generalized optimizers, and the morphisms between the generalized optimizers $d: A^{k}\rightarrow A^{k}$ and $d': B^{k}\rightarrow B^{k}$ are $\base$-morphisms $f: A\rightarrow B$ such that $f^{k} \circ d = d' \circ f^k: A^k \rightarrow B^k$. Note that morphisms only exist between generalized optimizers with the same dimension. The composition of morphisms in $\optimizer$ is the same as in $\base$. 
\end{definition}
Recall that $A^{k}$ and $f^{k}$ are respectively $A$ and $f$ tensored with themselves $k$ times. In the standard domain a generalized optimizer with dimension $k$ is a tuple $(\rl^n,d)$ where $d:\rl^{kn} \rightarrow \rl^{kn}$ is an optimizer (Definition \ref{optimizerdefinition}).
\begin{definition}\label{continuousoptimizationfunctordefinition}
Given a subcategory $\db$ of $\objective$, an \textbf{optimization functor over $\db$} is a functor $\db \rightarrow \optimizer$ that maps the objective $l: A \rightarrow X$ to a generalized optimizer over $(\base, X)$ with state space $A$.
\end{definition}
Optimization functors are generalizations of optimization schemes (Definition \ref{optimizationschemedefinition}) that map objectives to generalized optimizers. Explicitly, an optimization scheme $u$ that maps $l: \rl^n \rightarrow \rl$ to $u(l): \rl^{kn} \rightarrow \rl^{kn}$ defines an optimization functor in the standard domain.

The invariance properties of optimization functors are represented by the subcategory $\db \subseteq \objective$ out of which they are functorial. Concretely, consider the following categories:
\begin{itemize}
    \item $\objectivei$: The subcategory of $\objective$ in which morphisms are limited to invertible linear morphisms $l$ in $\base$.
    \item $\objectiveo$: The subcategory of $\objectivei$ in which the inverse of $l$ is $l^{\dagger}$.
\end{itemize}
In both the standard domain and $r$-polynomial domain, the morphisms in $\objectivei$ are linear maps defined by an invertible matrix and the morphisms in $\objectiveo$ are linear maps defined by an orthogonal matrix (matrix inverse is equal to matrix transpose). We will now generalize Proposition \ref{invarianceprop} by defining generalized gradient descent and momentum functors that are functorial out of $\objectiveo$ and a generalized Newton's method functor that is functorial out of $\objectivei$.
\begin{definition}
\textbf{Generalized gradient descent} sends the objective $l: A \rightarrow X$ to the generalized optimizer $ -R[l]_1: A \rightarrow A$ with dimension $1$.
\end{definition}
\begin{definition}
\textbf{Generalized momentum} sends the objective $l: A \rightarrow X$ to the generalized optimizer $\langle \pi_1, - \pi_1 - (R[l]_1 \circ \pi_0)\rangle: A^2 \rightarrow A^2$ with dimension $2$.
\end{definition}
Generalized momentum and generalized gradient descent have a very similar structure, with the major difference between the two being that generalized momentum uses a placeholder variable and generalized gradient descent does not. In the standard domain we have that $-R[l]_1(x) = -\nabla l(x)$ and $(\langle \pi_1, -\pi_1 - (R[l]_1 \circ \pi_0)\rangle)(x,y) = (y, -y - \nabla l (x))$, so generalized gradient descent and generalized momentum are equivalent to the gradient descent and momentum optimization schemes that we defined in Section \ref{continuousoptimizationalgorithm}. Similarly, in the $r$-polynomial domain generalized gradient descent maps $l_P:n \rightarrow 1$ to $-R[l_P]_1: n \rightarrow n$ and generalized momentum maps $l_P$ to $\langle \pi_1, -\pi_1 - (R[l_P]_1 \circ \pi_0) \rangle: n^2 \rightarrow n^2$ where:
\begin{gather*}
    \langle \pi_1, -\pi_1 - (R[l_P]_1 \circ \pi_0) \rangle (x, x') = 
    \left(x',
    -x' - \left(
    \frac{\partial l_P}{\partial x_1}(x),
    \cdots, 
    \frac{\partial l_P}{\partial x_n}(x)
    \right)\right)
\end{gather*}
Since Newton's method involves the computation of an inverse Hessian it is not immediately obvious how we can express it in terms of Cartesian reverse derivatives. However, by the inverse function theorem
we can rewrite the inverse Hessian as the Jacobian of the inverse gradient function, which makes this easier. That is:
\begin{gather}\label{jacobianhessian}
    %
    %
    %
    (\nabla^2 l)(x)^{-1} = 
    %
    J_{\nabla l}(x)^{-1} =
    %
    J_{(\nabla l)^{-1}}(\nabla l(x))
\end{gather}
where $J_{\nabla l}(x) = (\nabla^2 l)(x)$ is the Hessian of $l$ at $x$, $J_{(\nabla l)^{-1}}(\nabla l(x))$ is the Jacobian of the inverse gradient function evaluated at $\nabla l(x)$, and the second equality holds by the inverse function theorem. We can therefore generalize the Newton's method term $-\nabla^2 (l )^{-1}\nabla l$ as $-R[R[l]_1^{-1}] \circ \langle R[l]_1, R[l]_1 \rangle: X \rightarrow X$ and generalize Newton's method as follows:
\begin{definition}
\textbf{Generalized Newton's method} sends the objective $l: A \rightarrow X$ to the generalized optimizer $ -R[R[l]_1^{-1}] \circ \langle R[l]_1, R[l]_1\rangle: A \rightarrow A$ with dimension $1$.
\end{definition}
Equation \ref{jacobianhessian} implies that generalized Newton's method in the standard domain is equivalent to the Newton's method optimization scheme that we defined in Section \ref{continuousoptimizationalgorithm}. In the $r$-polynomial domain generalized Newton's Method maps the polynomial $l_P:n \rightarrow 1$ to $-R[R[l_P]_1^{-1}] \circ \langle R[l_P]_1, R[l_P]_1\rangle: n \rightarrow n$ where:
\begin{align*}
    -R[R[l_P]_1^{-1}] \circ \langle R[l_P]_1, R[l_P]_1\rangle (x) =  \\
    -\begin{pmatrix} 
    \frac{\partial (R[l_P]_{1}^{-1})_{1}}{\partial x_1}(R[l_P]_1(x))
    & \dots & 
    \frac{\partial (R[l_P]_{1}^{-1})_{1}}{\partial x_n}(R[l_P]_1(x)) \\
    \vdots & \ddots & \\
    \frac{\partial (R[l_P]_{1}^{-1})_{n}}{\partial x_1}(R[l_P]_1(x))
    &  &
    \frac{\partial (R[l_P]_{1}^{-1})_{n}}{\partial x_n}(R[l_P]_1(x))
    \end{pmatrix}^{T}
    \left(
    \frac{\partial l_P}{\partial x_1}(x),
    \cdots, 
    \frac{\partial l_P}{\partial x_n}(x)
    \right)
\end{align*}
Note that $(R[l_P]_{1}^{-1})_{i}$ is the $i$th projection of the inverse of the reverse derivative map. We now generalize Proposition \ref{invarianceprop}:
\begin{proposition}\label{generalizedinvarianceprop}
Generalized Newton's method is a functor from $\objectivei$ to $\optimizer$, whereas both generalized gradient descent and generalized momentum are functors from $\objectiveo$ to $\optimizer$.
(Proof in Appendix \ref{generalizedinvarianceproof})
\end{proposition}
Proposition \ref{generalizedinvarianceprop} implies that the invariance properties of our optimization functors mirror the invariance properties of their optimization scheme counterparts. Not only does Proposition \ref{generalizedinvarianceprop} directly imply Proposition \ref{invarianceprop}, but it also implies that the invariance properties that gradient descent, momentum, and Newton's method enjoy are not dependent on the underlying category over which they are defined.

\subsection{Generalized Optimization Flows}\label{generalizedoptimizationflow}


In Section \ref{standardoptimization} we demonstrated how we can derive continuous and discrete dynamical systems from an optimizer $d: \rl^{kn} \rightarrow \rl^{kn}$. In this section we extend this insight to generalized optimizers.

To do this, we define a morphism $s: X \rightarrow A^k$ whose Cartesian derivative is defined by a generalized optimizer $d: A^k \rightarrow A^k$. Since we can interpret morphisms in $\base[*, X]$ as either times $t$ or objective values $x$, the morphism $s: X \rightarrow A^k$ describes how the state of our dynamical system evolves in time.
Formally we can put this together in the following structure:
\begin{definition}
A \textbf{generalized optimization flow} over the optimization domain $\basedomain$ with \textbf{state space} $A \in \base$ and \textbf{dimension} $k \in \mathbb{N}$ is a tuple $(l, d, s, \tau)$ where $l: A \rightarrow X$ is an objective, $d: A^k \rightarrow A^k$ is a generalized optimizer, $s: X \rightarrow A^k$ is a morphism in $\base$ and $\tau$ is an interval in $\base[*, X]$ such that for $t \in \tau$ we have
%
$d \circ s \circ t =  D_1[s] \circ t: * \rightarrow A^{k}$.
\end{definition}
Intuitively, $l$ is an objective, $d$ is a generalized optimizer, and $s$ is the \textbf{state map} that maps times in $\tau$ to the system state such that $d \circ s: X \rightarrow A^{k}$ describes the Cartesian derivative of the state map $D_1[s]$.

In the standard domain we can define a generalized optimization flow $(l, d, s, \rl)$ from an optimizer $d:\rl^{kn} \rightarrow \rl^{kn}$ and an initial state $s_0 \in \rl^{kn}$ by defining a state map $s: \rl \rightarrow \rl^{kn}$ where $s(t) = s_0 + \int_{0}^{t} d(s(t')) dt'$. We can think of a state map in the standard domain as a simulation of Euler's method with infinitesimal $\alpha$:
\begin{gather*}
lim_{\alpha \rightarrow 0} s(t + \alpha)  =
lim_{\alpha \rightarrow 0}  s(t) + \alpha d(x_t)
\end{gather*}
%

\begin{definition}
A generalized optimization flow $(l, d, s, \tau)$ over the optimization domain $\basedomain$ is an \textbf{$n$-descending flow} if for any $t \in \tau$ and $k\leq n$ we have:
\begin{gather*}
    D_k[l \circ \pi_0 \circ s] \circ t \leq 0_X: * \rightarrow X
\end{gather*}
\end{definition}
Note that if $(l, d, s, \tau)$ is an $n$-descending flow and $l \circ \pi_0 \circ s: X \rightarrow X$ is $n$-smooth (Definition \ref{nsmooth}), then $l \circ \pi_0 \circ s$ must be monotonically decreasing in $t$ on $\tau$.

\begin{definition}
%
The generalized optimization flow $(l, d, s, \tau)$ over the optimization domain $\basedomain$ \textbf{converges} if for any $\delta > 0_X: * \rightarrow X$ there exists some $t \in \tau$ such that for any $t \leq t' \in \tau$ we have $-\delta \leq (l \circ \pi_0 \circ  s \circ t') - (l \circ \pi_0 \circ  s \circ t) \leq \delta$

\end{definition}
In the standard domain this reduces to a familiar definition of convergence that is similar to what \cite{andersongradientflow} uses: a flow converges if there exists a time $t$ after which the value of the objective $l$ does not change by more than an arbitrarily small amount.

Now suppose $(l, d, s, \tau)$ is an $n$-descending flow, $l \circ \pi_0 \circ s: X \rightarrow X$ is $n$-smooth and $l$ is bounded below (Definition \ref{boundedbelow}). Since $l \circ \pi_0 \circ s$ must decrease monotonically in $t$ it must be that $(l, d, s, \tau)$ converges. In the next section we give examples of optimization flows defined by the generalized gradient that satisfy these conditions.

\subsubsection{Generalized Gradient Flows}

\begin{definition}
A \textbf{generalized gradient flow} is a generalized optimization flow of the form $(l, -R[l]_1, s, \tau)$.
\end{definition}
Given a smooth objective $l:\rl^n \rightarrow \rl$ an example generalized gradient flow in the standard domain is $(l, -\nabla l, s, \rl)$ where $s(t) = s_0 + \int_{0}^{t} -\nabla l(s(t')) dt'$ for some $s_0 \in \rl^n$. One of the most useful properties of a generalized gradient flow is that we can write its Cartesian derivative with an inner product-like structure:
\begin{proposition}\label{innerproductproposition}
Given a choice of time $t\in\tau$ and a generalized gradient flow $(l, -R[l]_1, s, \tau)$ we can write the following:
\begin{gather*}
    D_1[l \circ \pi_0 \circ s] \circ t = -R[l]_{s_t}^{\dagger}
    \circ 
    R[l]_{s_t} \circ
    1_X: * \rightarrow X
\end{gather*}
where $R[l]_{s_t} = R[l] \circ \langle s \circ t \circ !_{X}, id_X \rangle: X \rightarrow A$.
(Proof in Appendix \ref{innerproductproof})
\end{proposition}
Intuitively, $s \circ t: * \rightarrow A$ is the state at time $t$ and $R[l]_{s_t} \circ 1_X: * \rightarrow A$ is the value of the generalized gradient of $l$ at time $t$. 
To understand the importance of this result consider the following definition:
\begin{definition}
$\basedomain$ \textbf{supports generalized gradient-based optimization} when any generalized gradient flow over $\basedomain$ is a $1$-descending flow.
\end{definition}
Intuitively, an optimization domain supports generalized gradient-based optimization if loss decreases in the direction of the gradient. Proposition \ref{innerproductproposition} is important because it helps us identify the optimization domains for which this holds. For example, Proposition \ref{innerproductproposition} implies that both the standard domain and any $r$-polynomial domain support generalized gradient-based optimization:
\begin{itemize}
    \item In the standard domain we have that:
    \begin{gather*}
        -R[l]_{s_t}^{\dagger}
        \circ 
        R[l]_{s_t} \circ 
        1_{\rl} =
        -\nabla l(s(t))^{T} \nabla l(s(t))  = 
        -\|\nabla l(s(t))\|^2
    \end{gather*}
    which must be non-positive by the definition of a norm. As a result, any generalized gradient flow $(l, -R[l], s, \tau)$ in the standard domain converges if $l$ is bounded below.
    \item In the $r$-polynomial domain we have that:
    \begin{gather*}
        -R[l_P]_{s_t}^{\dagger}
        \circ 
        R[l_P]_{s_t} \circ 
        1_1
        = 
        %
        %
        -\sum_{i=1}^n \frac{\partial l_P}{\partial x_i}(s_{t}) 
        R[l_P]_{s_t}(1_1)_i
        =
        %
        -\sum_{i=1}^n \frac{\partial l_P}{\partial x_i}(s_{t}) 
        \frac{\partial l_P}{\partial x_i}(s_t)
    \end{gather*}
    %
    which must be non-positive since in an ordered ring no negative element is a square. If $r$ is a dense subring of a real-closed field then any generalized gradient flow $(l, -R[l], s, \tau)$ in the $r$-polynomial domain converges if $l$ is bounded below  (since $l\circ \pi_0 \circ s: X \rightarrow X$ must be $1$-smooth, see Section \ref{generalizedgradientderivativesection}).
\end{itemize}

\section{Example and Experiment}\label{experiments}

We start this section with a demonstration of the structure and behavior of an example optimization flow. We then build on this example to define an algorithm for finding integer minima of multivariate polynomials. We demonstrate that this algorithm consistently outperforms random search.

\subsection{Illustrative Example - Integer Polynomial State Map}

Suppose $l_P: 1 \rightarrow 1$ is an objective and $u$ is an optimization functor in the integer polynomial domain $\polyz$. Given a choice of integer $x_0 \in \mathbb{Z}$, we can follow the pattern laid out in Section \ref{standardoptimization} and form a discrete dynamical system $x_{t+1} = x_t + u(l_P)(x)$. 

Now suppose that for some $\tau = 1,2,\cdots,m$ we would like to construct an optimization flow $(l_P, u(l_P), s, \tau)$ that traces out the values of this dynamical system. The state map $s$ must be an integer polynomial that satisfies two properties:
\begin{enumerate}
    \item The integer polynomial $s$ intersects the discrete dynamical system at each $t \in \tau$:
    \begin{gather*}
        s(t+1) = s(t) + u(l_P)(s(t))
    \end{gather*}
    \item By the definition of an optimization flow it must be that $u(l_P)$ defines the derivative of $s$. That is, for  $t \in \tau$:
    \begin{gather*}
        u(l_P)(s(t)) = \frac{\partial s}{\partial t}(t)
    \end{gather*}
\end{enumerate}
By Proposition \ref{innerproductproposition} we expect that $s$ will move towards the minima of $l_P$ at each step.

There may be no, some, or an infinite number of integer polynomials $s(t) = p_1 t + p_2 t^2 + \cdots + p_n t^n$ that satisfy these conditions. For example, consider the simple case in which $l_P(x) =  a x^2 + b$ and $u(l_P) = -R[l_P] = -\frac{\partial l_P}{\partial x}$. In this case the condition $u(l_P)(s(t)) - \frac{\partial s}{\partial t}(t) = 0$ becomes:

\begin{align*}
    u(l_P)(s(t)) - \frac{\partial s}{\partial t}(t) = 
    \\
    -\frac{\partial l_P}{\partial x}(s(t)) - \frac{\partial s}{\partial t}(t) =
    \\
    -2 a s(t) - \frac{\partial s}{\partial t}(t) =
    \\
    -\left(2 a p_0 + 2 a p_1 t + 2 a p_2 t^2 + \cdots + 2 a p_n t^n\right) - \left(
    p_1 + 2 p_2 t + 3 p_3 t^2 + \cdots + n p_n t^{n-1}\right) =
    \\
    (2 a p_0 + p_1) + (2 a p_1 + 2 p_2) t + (2 a p_2 + 3 p_3) t^2 + \cdots + (2 a p_n + n p_n) t^{n-1} + 2 p_n t^n  = 0
\end{align*}
and $s(t+1) - s(t) - u(l_P)(s(t)) = 0$ becomes:
\begin{align*}
    s(t+1) - s(t) - u(l_P)(s(t)) = 
    \\
    s(t+1) - s(t) + \frac{\partial l_P}{\partial x}(s(t)) = 
    \\
    s(t+1) - s(t) + 2 a s(t) = 
    \\
    s(t+1) + (2a - 1)s(t) = 
    \\
    \left(p_0 + p_1 (t+1) + \cdots + p_n (t + 1)^n
    \right) + \left(
    (2a - 1) p_0 + (2a - 1) p_1 t + \cdots + (2a - 1) p_n t^n\right)= 
    \\
    \sum_{i=0}^n p_i ((t+1)^{i} + (2a - 1) t^{i}) = 0
\end{align*}

\begin{figure}[ht]\begin{center}
\includegraphics[width=12cm,height=6cm]{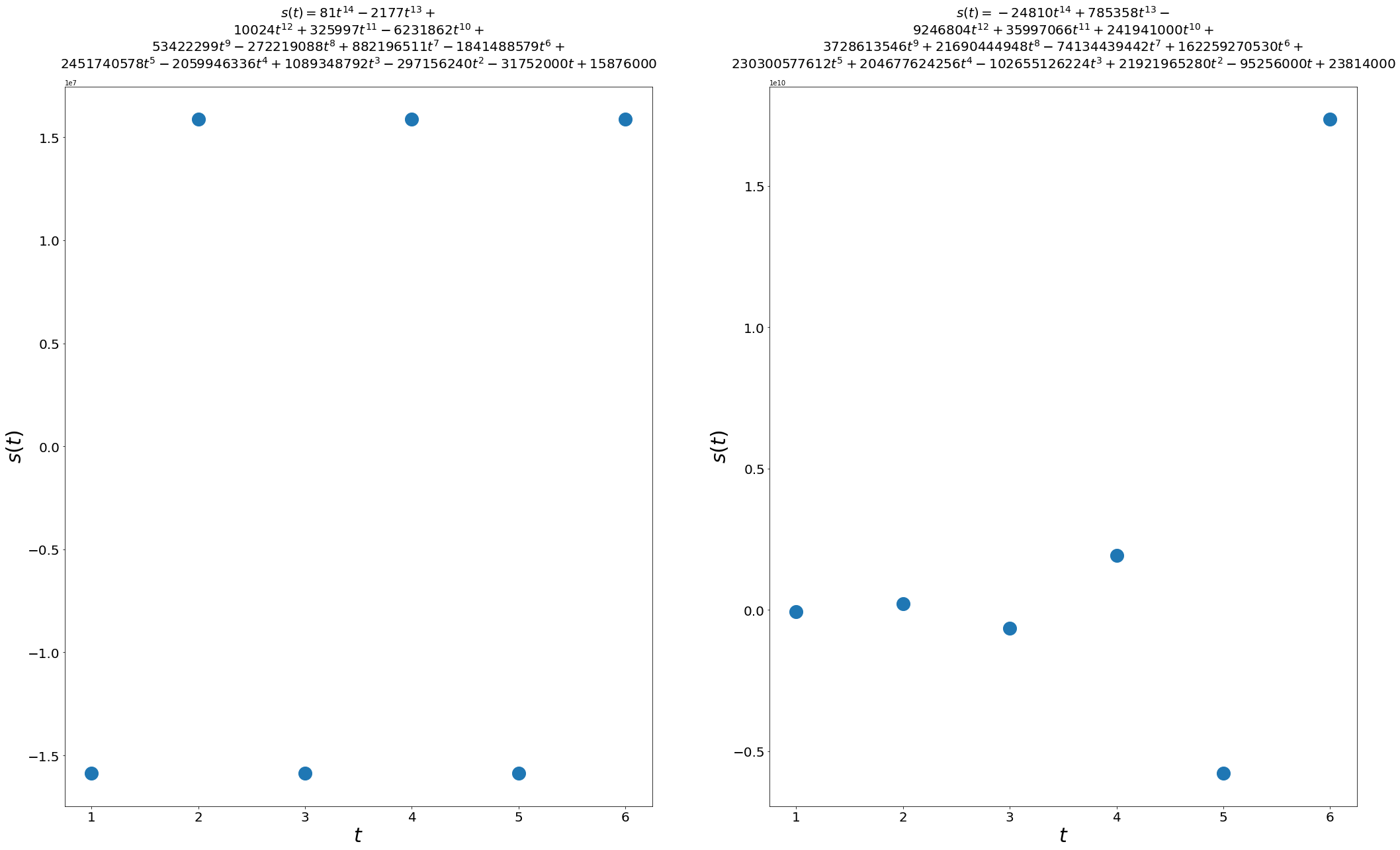}
\caption{
Example integer polynomial state maps $s(t) = p_1 t + p_2 t^2 + \cdots + p_n t^n$ for $\tau=1,2,\cdots,6$. If $l_P = x^2$ (left) and $l_P = 2x^2 - 1$ (right) then $(l_P, -R[l_P]_1, s, \tau)$ forms a gradient flow.  }
\label{trajectories}
\end{center}
\end{figure}

Evaluated at each $t=1,2,\cdots,m$ this forms a linear Diophantine system with $2m+1$ unique equations. There are therefore infinitely many degree $2m+2$ polynomials $s$ that satisfy these equations. We show two examples in Figure \ref{trajectories}. As we would expect from Proposition \ref{innerproductproposition}, we see that each step the dynamical system takes is in the right direction.

\subsection{Experiment - Integer Gradient Descent}


Although Figure \ref{trajectories} shows that each step the dynamical system takes is in the right direction, these steps are too large to minimize the function, which can cause $s$ to oscillate or diverge. In the standard domain we would mitigate this problem by choosing a smaller step size (aka learning rate) $\alpha$ for the dynamical system $x_{t+1} = x_t - \alpha \frac{\partial l_P}{\partial x}(x_t)$, but this is not possible in the integer polynomial domain since there are no $\alpha$ between $0$ and $1$. Instead, we can simply modify the dynamical system to instead take steps of size $1$ in the direction of the negative gradient.


We can assess how well this method performs at finding integer minima of arbitrary multivariate polynomials by testing it on randomly generated polynomials. In Table \ref{polynomialtable} we show that this method consistently outperforms random search at finding minima of polynomials that can be written as a sum of squared terms (which are guaranteed to have both a global minimum and a global integer minimum).


\begin{table}
    \centering
    \begin{tabular}{|c|c|c|c|c|}
        \hline
        \multicolumn{1}{|p{2cm}|}{\centering Number of Steps ($N$)}
        &
        \multicolumn{1}{|p{5.5cm}|}{\centering Frequency that Integer Gradient Descent is Better than Random Search}
        \\
        \hline
        \hline
        5 &  0.740 ($\pm$ 0.019) \\
        \hline
        10   & 0.763 ($\pm$ 0.019) \\
        \hline
        50   & 0.769 ($\pm$ 0.020) \\
        \hline
        100   & 0.807 ($\pm$ 0.011) \\
        \hline
    \end{tabular}
    \caption{In this experiment we randomly generate polynomials $l_P$ as sums of squared terms and use both integer gradient descent and random search to minimize $l_P$. In integer gradient descent we sample a point and take $N$ steps in the gradient direction. In random search we sample $N$ points (log-uniformly) and choose the best one. We then compute the frequency with which gradient descent finds a better $x$ than random search. The number of variables, terms, and coefficient values are sampled from $[1,10]$. Means and standard errors from $100$ experiments of $10$ polynomials each are shown. The code to run these experiments is on Github at tinyurl.com/ku3pjz56.
    }
\label{polynomialtable}
\end{table}

\section{Discussion and Future Work}\label{discussion}

In recent years researchers have begun to study categorical generalizations of machine learning. This research has proceeded on many fronts:
\cite{fritz2020synthetic} and \cite{cho2019disintegration} introduce synthetic perspectives on probability theory, \cite{fong2019backprop} introduce a functorial perspective on the backpropagation algorithm and  \cite{autodiff2018, blute2009cartesian, cruttwell2021categorical} explore categorical formulations of automatic differentiation. Some authors have also begun to explore categorical generalizations of classical machine
learning techniques. \cite{cho2019disintegration} introduce a generalized perspective on Bayesian updating and \cite{Wilson_2021} introduce a generalized perspective on gradient descent that can be used to learn Boolean circuits.

Despite this progress, there has been relatively little research on the properties of these generalized algorithms. That is, although categorical machine learning has started to gain traction, categorical learning theory is still far behind. In this paper we aim to reduce this gap by exploring the properties of optimizers generalized over other categories.

However, there is still much to do. For example, although we identify the properties that a generalized optimization flow must possess in order to converge, our construction does not distinguish between flows that converge to minima or to arbitrary points. Furthermore, there are many variations of gradient-based optimizers that our formulation does not capture, such as stochastic optimizers like stochastic gradient descent and mini-batch gradient descent.

Another potential future direction for this work is to explore generalizations of constrained optimization. Optimization algorithms like gradient descent and Newton's method can be adapted to solve constrained optimization problems, and we may be able to do the same for their generalized analogs. This may enable us to adapt the technique for minimizing integer polynomials that we introduce in Section \ref{experiments} towards solving integer programs, which is an NP-hard problem with an enormous number of practical applications.


We believe that this line of research will accelerate future machine learning research by helping researchers better understand the foundational components of the algorithms that they use. This generalized perspective may also help us better understand the domains over which different algorithms will be successful.


\bibliographystyle{plainnat}
\bibliography{main.bib}
\appendix


\section{Appendix - Proofs}\label{9-appendix}

\subsection{Proof of Proposition \ref{invarianceprop}}\label{invarianceproof}
\begin{proof}

First, we will show that the Newton's method optimizer scheme $NEW(l)(x) = -(\nabla^2 l (x))^{-1}\nabla l(x)$ is invariant to invertible linear transformations. Consider any function of the form $f(x) = Ax$ where $A$ is invertible. We have:

\begin{align*}
NEW(l \circ f)(x) = \\
-(\nabla^{2} (l \circ f)(x))^{-1} \nabla (l \circ f)(x) =\\
-A^{-1}(\nabla^{2}l(Ax))^{-1} A^{-T}A^{T}\nabla l(Ax) = \\
-A^{-1}(\nabla^{2}l(Ax))^{-1}\nabla l(Ax) = \\
-f^{-1}((\nabla^{2}l(f(x)))^{-1} \nabla l(f(x))) =\\
f^{-1}(NEW(l)(f(x)))
\end{align*}
Next, we will show that the gradient descent optimizer scheme $GRAD(l)(x) = \nabla l (x)$ is invariant to orthogonal linear transformations, but not to linear transformations in general. Consider any function of the form $f(x) = Ax$ where $A$ is an orthogonal matrix. Then the following holds only when $A^{T} = A^{-1}$:
\begin{align*}
GRAD(l \circ f)(x) =  \\
-\nabla (l \circ f)(x) = \\
-A^{T} (\nabla l (Ax)) = \\
-A^{-1} (\nabla l (Ax)) = \\
-f^{-1}(GRAD(l)(f(x)))
\end{align*}
Next, we will show that the momentum optimizer scheme $MOM(l)(x, y) = (y, y + \nabla l(x))$ is also invariant to orthogonal linear transformations, but not to linear transformations in general. Consider any function of the form $f(x) = Ax$ where $A$ is an orthogonal matrix. Then the following holds only when $A^{T} = A^{-1}$:
\begin{align*}
MOM(l \circ f)(x, y)_x = y = A^{T}Ay = f^{-1}(MOM(l)(f(x), f(y)))_x
\\
MOM(l \circ f)(x, y)_y =
- y - \nabla (l \circ f)(x)) =
- A^{-1}Ay - A^{T} \nabla l(Ax)) =
f^{-1}(MOM(l)(f(x), f(y)))_y
\end{align*}

\end{proof}

\subsection{Proof of Proposition \ref{eulersmethodprop}}\label{eulersmethodproof}
\begin{proof}
Consider starting at some point $x_0 \in \mathbb{R}^{kn}$ and repeatedly taking Euler steps $x_{t+\alpha} = x_t + \alpha u(l)(x_t)$. Now suppose instead that we start at the point $y_0 = f_{k}^{-1}x_0$ and take Euler steps $y_{t+\alpha} = y_t + \alpha u(l \circ f)(y_t)$. 

We will prove by induction that $y_{t+\alpha} = f_{k}^{-1}(x_{t+\alpha})$, and therefore the two sequences converge at the same rate. The base case holds by definition and by induction we can see that:
\begin{align*}
y_{t+\alpha} =
y_t + \alpha u(l \circ f)(y_t) =
%
f_{k}^{-1}(x_t) + \alpha f_{k}^{-1}(u(l)(x_t)) =
%
f_{k}^{-1}(x_{t+\alpha})
\end{align*}
\end{proof}


\subsection{Proof of Proposition \ref{generalizedinvarianceprop}}\label{generalizedinvarianceproof}
\begin{proof}
Since generalized gradient descent, generalized momentum and generalized Newton's method all act as the identity on morphisms, we simply need to show that each functor maps a morphism in its source category to a morphism in its target category.

First we show that generalized Newton's method $NEW(l) = R[R[l]_1^{-1}] \circ \langle R[l]_1, R[l]_1 \rangle$ is a functor out of $\objectivei$. Given an objective $l: A \rightarrow X$ and an  invertible linear map $f: B \rightarrow A$ we have:
\begin{gather*}
    f \circ NEW(l \circ f ): B \rightarrow A
\end{gather*}
\begin{align*}
f \circ NEW(l \circ f )= \\ 
%
-f \circ R[R[l \circ f]_1^{-1}] \circ \langle R[l \circ f]_1, R[l \circ f]_1 \rangle
=^{*}  \\
%
-f \circ f^{-1} \circ R[R[l]_1^{-1}]
\circ 
(f^{-\dagger} \times f^{-\dagger})
\circ 
\langle
f^{\dagger} \circ R[l]_1 \circ f, f^{\dagger} \circ R[l]_1 \circ f
\rangle
=  \\
%
-R[R[l]_1^{-1}]
\circ 
\langle
f^{-\dagger} \circ f^{\dagger} \circ R[l]_1 \circ f,
f^{-\dagger} \circ f^{\dagger} \circ R[l]_1 \circ f
\rangle
=  \\
%
-R[R[l]_1^{-1}]
\circ 
\langle R[l]_1, R[l]_1 \rangle
\circ
f
= \\
%
NEW(l) \circ f
\end{align*}
where $*$ holds by:
\begin{gather*}
    R[R[l \circ f]_{1}^{-1}]: B \times B \rightarrow B
\end{gather*}
\begin{align*}
R[R[l \circ f]_{1}^{-1}]  =^{**}  \\ 
%
R[f^{-1}  \circ R[l]_1^{-1} \circ f^{-\dagger}] = \\
%
R[f^{-\dagger}] 
\circ 
(id_B \times R[f^{-1}  \circ R[l]_1^{-1}]) 
\circ 
\langle \pi_0, \langle f^{-\dagger} \circ \pi_0, \pi_1 \rangle \rangle 
= \\
%
f^{-1} 
\circ 
R[f^{-1}  \circ R[l]_1^{-1}] 
\circ 
\langle f^{-\dagger} \circ \pi_0, \pi_1 \rangle 
= \\
%
f^{-1} 
\circ 
R[R[l]_1^{-1}] 
\circ
(id_A \times R[f^{-1}]) 
\circ
\langle \pi_0,  \langle R[l]_1^{-1} \circ \pi_0, \pi_1 \rangle \rangle 
\circ 
\langle f^{-\dagger} \circ \pi_0, \pi_1 \rangle 
= \\
%
f^{-1} 
\circ 
R[R[l]_1^{-1}] 
\circ
(id_A \times f^{-\dagger}) 
\circ
\langle f^{-\dagger} \circ \pi_0, \pi_1 \rangle 
= \\
%
f^{-1} 
\circ 
R[R[l]_1^{-1}] 
\circ
(f^{-\dagger} \times f^{-\dagger})  
\end{align*}
and where $**$ holds by:
\begin{gather*}
    R[l \circ f]_{1}^{-1}: B\rightarrow B 
\end{gather*}
\begin{align*}
    R[l \circ f]_{1}^{-1}  =
    f^{-1} \circ R[l]_1^{-1} \circ R[f^{-1}] \circ ( 1_A \times id_{B})   =
    f^{-1} \circ R[l]_1^{-1} \circ f^{-\dagger}
\end{align*}

Next we show that generalized gradient descent $GRAD(l)  = (1, A, R[l]_1)$ is a functor out of $\objectiveo$. Given an objective $l: A \rightarrow X$ and an invertible linear map $f: B \rightarrow A$ where $f \circ f^{\dagger} = id_A$ and $f^{\dagger} \circ f = id_B$ we have:
\begin{gather*}
    f \circ GRAD(l\circ f): B \rightarrow A
\end{gather*}
\begin{align*}
f \circ GRAD(l\circ f)  =\\
-f \circ R[l \circ f]_1  =\\
%
-f \circ R[l \circ f] \circ \langle id_{B}, 1_{BX} \rangle=\\
%
-f\circ R[f] \circ (id_{B} \times R[l]_1) \circ \langle id_B, f \rangle =\\
%
-f\circ f^{\dagger} \circ \pi_1 \circ (id_{B} \times R[l]_1) \circ \langle id_B, f \rangle =\\
%
-\pi_1 \circ   (id_B \times R[l]_1) \circ \langle id_B, f \rangle =\\
%
-R[l]_1 \circ f =\\
%
GRAD(l) \circ f
\end{align*}
Next we show that generalized momentum $MOM(l) = (1, A, \langle \pi_1, \pi_1 + (R[l]_1 \circ \pi_0)\rangle)$ is a functor out of $\objectiveo$. Given an objective $l: A \rightarrow X$ and an invertible linear map $f: B \rightarrow A$ where $f \circ f^{\dagger} = id_A$ and $f^{\dagger} \circ f = id_B$ we have:
\begin{gather*}
    f \circ MOM(l\circ f): B^2 \rightarrow A^2
\end{gather*}
\begin{align*}
f^{2} \circ MOM(l \circ f) = \\
(f \times f) \circ MOM(l \circ f) = \\
(f \times f) \circ \langle \pi_1, - \pi_1 - (R[l \circ f]_1 \circ \pi_0)\rangle = \\
%
=\langle f \circ \pi_1, f \circ (-\pi_1 - (R[l \circ f]_1 \circ \pi_0)) \rangle = \\
%
\langle f \circ \pi_1, -f \circ \pi_1 - (f \circ  R[l \circ f]_1 \circ \pi_0) \rangle = \\
%
\langle f \circ \pi_1, -f \circ \pi_1 - (R[l]_1 \circ f \circ \pi_0) \rangle = \\
%
\langle \pi_1, -\pi_1 - (R[l]_1 \circ \pi_0)\rangle \circ (f \times f) =
\\
MOM(l) \circ (f \times f) =
\\
MOM(l) \circ f^{2}
\end{align*}

\end{proof}


\subsection{Proof of Proposition \ref{innerproductproposition}}\label{innerproductproof}
\begin{proof}
For $D_1[l \circ \pi_0 \circ s] \circ t : 1 \rightarrow X$ we have that:
\begin{align*}
    D_1[l \circ \pi_0 \circ s] \circ t = \\
    %
    D_1[l \circ s] \circ t = \\
    %
    D[l \circ s] \circ \langle t, 1_{X}\rangle = \\
    %
    D[l] \circ \langle s \circ \pi_0, D[s]\rangle \circ \langle t, 1_{X}\rangle  = \\ 
    %
    D[l] \circ \langle s, D[s] \circ \langle id_X, 1_{X}\rangle \rangle \circ t  = \\ 
    %
    D[l] \circ \langle s, d \circ s \rangle \circ t = \\
    %
   D[l] \circ \langle s,  -R[l] \circ \langle id_{A}, 1_{AX} \rangle \circ s \rangle  \circ t 
    = \\
    %
   -D[l] \circ \langle s,  R[l] \circ \langle s, 1_{X}\rangle \rangle  \circ t 
    = \\
    %
   -\pi_1
   \circ
   R[R[l]]
   \circ
   (\langle id_A, 1_{AX}  \rangle \times id_A)
   \circ
   \langle s,  R[l] \circ \langle s, 1_{X}\rangle \rangle  \circ t 
    = \\
    %
   -\pi_1
   \circ
   R[R[l]]
   \circ
   \langle \langle s, 1_{X} \rangle ,  R[l] \circ \langle s, 1_{X}\rangle \rangle  \circ t 
    = \\
    %
    -\pi_1
    \circ
    R[R[l]]
    \circ
    \langle \langle s \circ t, 1_{X}\rangle,  R[l] \circ \langle s, 1_{X}\rangle \rangle 
    = \\
    %
    -\pi_1
    \circ
    R[R[l]]
    \circ
    (\langle s \circ t, 1_{X}\rangle \times id_A)
    \circ 
    R[l] \circ \langle s \circ t, 1_{X}\rangle
    = \\
    %
    -\pi_1
    \circ
    R[R[l]]
    \circ
    (\langle s \circ t, 1_{X}\rangle \times id_A)
    \circ 
    R[l]_{s_t} \circ 
    1_{X}
    = \\
    %
    %
    -R[R[l] \circ \langle s \circ t \circ !_{X}, id_X \rangle ]
    \circ 
    \langle 1_{X}, R[l]_{s_t}\rangle
    \circ 
    1_{X} = \\
    %
    -(R[l] \circ \langle s \circ t \circ !_{X}, id_X \rangle )^{\dagger}
    \circ 
    \pi_1
    \circ 
    \langle 1_{X}, R[l]_{s_t}\rangle
    \circ 
    1_{X}
    = \\
    %
    -(R[l] \circ \langle s \circ t \circ !_{X}, id_X \rangle )^{\dagger}
    \circ 
    R[l]_{s_t}
    \circ 
    1_{X}
    = \\
    %
    -R[l]_{s_t}^{\dagger}
    \circ 
    R[l]_{s_t} \circ 
    1_{X}
\end{align*}
\end{proof}

\end{document}